\documentclass{matstos}
\usepackage[OT4,plmath]{polski}
\usepackage[cp1250]{inputenc}

\received{}
\lastrevision{}
\firstpage{35}    
\volume{13}
\fasc{54}
\years{2011}
\pages{35--46}

\logo{}{}{}
\received{}

\title[Applications of topology in computer algorithms]{Applications of topology in computer algorithms}

\author[Rastislav Telg\'arsky]{Rastislav Telgarsky }
\thanks{This paper is based on the invited lecture at International Conference on Topology and Applications held in August 23--27, 1999, at Kanagawa University in Yokohama, Japan}
\affiliation{Department of Mathematics}
\address{Central New Mexico Community College}
\address{525 Buena Vista Dr. SE}
\address{Albuquerque, NM 87106, U.S.A.}
\email{rastislav@telgarsky.com}


\begin{document}
\setcounter{page}{35}
\begin{abstract}
The aim of this paper is to discuss some applications of general topology in computer algorithms including modeling and simulation, and also in computer graphics and image processing. While the progress in these areas heavily depends on advances in computing hardware, the major intellectual achievements are the algorithms. The applications of general topology in other branches of mathematics are not discussed, since they are not applications of mathematics outside of mathematics.
\end{abstract}

\section{Looking back}
After Ph.D. in Mathematics in 1970, I was 18 years in academic employments. Then I accepted a position related to algorithm development in a small company (about 25 employees). This kind of job, an industrial mathematician, was not recognized by the American Mathematical Society. The tenure avenues did not exist, and I had to cope every day with job insecurity. In order to survive and support my family I had to rebuild my education from my junior years up. In the evenings I studied ``engineering mathematics'' which I needed right away for the projects. The book that helped me the most was Numerical Recipes in C (Press et al.). The most helpful journals were IEEE Transactions on Pattern Analysis and Machine Intelligence (PAMI) and Computer Graphics and Image Processing (CGIP). My job position was described as senior mathematician, but other names could be software mathematician, mathematics engineer or algorithm engineer. After 8 years, in 1996, I joined a fast growing company, which was after few years bought by the Boeing Company. After being a~mathematician, a software engineer and an image processing scientist I was nominated to become Boeing Associate Technical Fellow. Only few years later, in 2007, I was included in 10\% layoff due to downsizing of the research and development base.

The main difference between the work in academia and industry is that in academia I do all research activities according to my choice and its limitations, while in industry, I am told what to do and I must do it as soon as possible. Thus, I must put away my desire to think about my questions or my study of a useful paper, because I must focus to the topics within the framework of the project. Usually, I was a member of an engineering team where members performed assigned and coordinated tasks. Obviously, in every industry organizations there are management and service positions. Due to circumstances, I landed in service positions, and the management positions were taken by electrical, mechanical, and software engineers, most of them with doctoral degrees.

\section{Mathematics in crisis.}
There is an on-going crisis in set-theoretic topology, in mathematics, and in scientific research in general. Two main factors are: the large number of candidates for every position, from 10:1 to over 1000:1, and insufficient funds for research and development.

In the U.S.A., general topology (set-theoretic topology) is not included in the Ivy League education. There are no set-theoretic topologists at Harvard, Yale and Princeton, while the differential and algebraic topologists still enjoy the academic acceptance. A number of Internet sites shows that the people attending the conferences in differential and algebraic topology are from the Ivy League Universities. Some journals, such as Annals of Mathematics and Pacific Journal of Mathematics, do not even publish papers in general topology. Many topologists, set-theorists, logicians, and model theorists of my generation went to computer science departments or to industry. General topology became a rare academic discipline and isolated from applications of mathematics.

Some engineering departments, such as electrical engineering, prefer that their faculty members teach the mathematics classes rather than letting the mathematics department delegate mathematicians to teach them. The engineers want the classes to be taught in a useful way, and not in a rigorous way mathematicians teach students of mathematics. As the result, many mathematics departments are in trouble.

Mathematics is usually applied by engineers or industrial scientists, and not by mathematicians. In other words, an applied mathematician does not apply mathematics; the engineers or scientists do it. In consequence, mathematicians are losing inspiration by problems of the real world. Many of them create an imaginary world, which ignores real-world problems. Also, they ignore the urgency of solving the real-world problems. Engineers took from us (mathematicians) the best and the most vital parts of mathematics! Should not we teach at least some students of mathematics how to be creative and useful in non-academic jobs? Should mathematicians accept the situation that most leading roles in industrial science research are by people having a major degree in engineering, physics, chemistry or biology?

Many scientists believe that their achievements are valid forever. They believe that science papers preserved on paper or recorded on other media will be available forever and that research will continue with references to previous results. However, scientific information is more like a genetic material: some survives and some dies. A lot of scientific information is becoming easy to find through the Internet. However, reaching the sources usually ends in an information flood. When I searched for ``scintillation'', the first 10 of about 1000 matches returned a web site with references to 300 published papers dealing with scintillation. Narrowing the search can be very difficult and time consuming!

Availability of research funds is a dream. What is the cost of writing a paper for publication, using library services/access, developing an algorithm for a client, building a device for a delivery, or even a patent preparation? How much do we need to know about the struggle for the funding of re\-search and development? Tenure positions create the class of people (including mathematicians) who ignore the needs of applied sciences. Also, there is the ever-presence of science politics. Which mathematics is important? Who gets the grant? Who gets the tenure? Who will stay and who will be pushed away? Who will teach 6 and who 16 hours per week? For how many people do you publish your paper? Will it appear in a low impact or in a high impact professional journal? How many referees will it pass through until it will get published? Is whom you know as important as the content in getting your paper published? Some people will minify and some will magnify your achievements. Who can be trusted?

I left the American Mathematical Society (AMS) in 1998 after 30 years of membership because it did not support my work as a mathematician in industry. In contrast, publications of Society of Photo-Optical Instrumentation Engineers (SPIE) were extremely helpful. For few years I held membership in the Society for Industrial and Applied Mathematics (SIAM), but I found it too narrowly focused, and not giving the slightest hint what a mathematician can or should do in an industrial employment. For few years I was a member of IEEE Computer Society, but only a segment of its activities was useful. I joined the Association for Computing Machinery (ACM) in 1986, but in few years, it was clear that the membership is irrelevant to my work. The world is changing very fast! Probably, a partial membership in several professional societies would give me an edge to become a mathematician in industry. The continuation of my employment was supported by written codes and by verified and approved/corrected codes of others.

\section{Limits of mathematical modeling}
How closely can mathematics describe reality? Does it only provide tools for an approximation of reality? We can think about infinite sets, but there are finite sets whose calculation would take more than 1000 years, which exceeds our lifetime. Nobody can make an arbitrarily small circle using a copper wire. Time slices cannot be arbitrarily small, and frequencies cannot be arbitrarily large. It turns out that every physical quantity has a limited range, and every natural law has a limited scope. Each real world sensor has a range; a sensor cannot detect arbitrarily small or arbitrarily large quantities.

Mathematics describes only some aspects of the real world. If mathematics approximates the laws of the real world, then which of our knowledge is exact and final, and which is always approximate and temporary? Are the representations of mathematical objects in computer circuits and in manufacturing them the final proofs of their consistency?

An arbitrary homeomorphism or even a contraction mapping does not make sense from a manufacturing point of view, because nobody can make a pencil of an arbitrary small size, in particular, of the size of a single atom. Our imagination of an arbitrarily small pencil does not match the physical reality. If we can manufacture an object A, and if A is homeomorphic to B, it is not necessarily possible to manufacture B. In topological terms, for example, if A is an iron ball, and B is a solid Alexander's Horned Sphere, then B would require infinitely many consecutive manufacturing steps with doubling complexity. It is necessary to accept the usage of both models of reality: the continuous models and the discrete (actually: finite) models. However, we are accepting mutually exclusive (and contradictory) views of the same reality, and we call them complementary.

Computers have become mathematical machines to compete with mathematicians in numerical and symbolic calculations. Also, the drawing as a way of describing an object (a construction) is being replaced by CAD computer graphics. However, the computer requires somewhat different mathematics to ``think''. Instead of arbitrarily small or large quantities we deal with small or large constants, respectively. For example, in numerical differentiation, the delta does not go to 0, but must stay above a small constant which depends on software and hardware configurations. The differential equations turn to difference equations, functions become vectors, operators become matrices, etc. The models of fluids and gases as continua are losing validity at a molecular level. Perfect white noise is a mathematical idealization and it only approximates the real white noise (or vice versa). Fractals with infinitely many levels do not exist in nature. The ferns in nature have only a few levels of self-similarity. The ferns (living plants) are the real things, not some approximations of real things.

\section{Topological notions and related algorithms}
The set-theoretic topology has been around for more than one hundred years, but no direct applications have been found. Let us see: if we accept the standard definition of topology as a family of open sets, or, as a study of invariants of homeomorphisms, then there are no applications of topology. In the real world the locations, the distances and sizes matter, thus it is the geometry and the metric topology that rule the world, although, sometimes with the help of statistics.

The boundary (frontier) of a set is defined via the closure operator. This concept is inconsistent with the physical reality. For example, in two layers of metals have two boundaries inside, corresponding to each layer. So, the boundary of one layer is disjoint with the boundary of its complements. Moreover, some physical sets do not have identifiable boundaries, for example, clouds, fogs, rains, etc.

The vital notions of (metric) topology such as connected sets, neighborhoods, boundaries (edges), surfaces, cluster points, distances, thinning, thickening, homotopy, and shape are used in various applications including image processing and pattern recognition. However, continuous functions (including homeomorphisms) and compact sets are useless in the same applications, where every set must be a finite list; every function is a vector, etc.).


The notions of classical analytic geometry such as lines, triangles, polygons, circles, ellipses, etc. lead to computational geometry and computer graphics. However, computational geometry is a digital (that is, discrete and finite) geometry. For example, the Bresenham line drawing algorithm creates line segments made of a small finite number of pixels. Ray tracing and ray casting methods are based on analytic solutions to intersection of lines and solids.

Although the origin of fractals was in topology (curves, continua, singular sets, etc.), there is no branch of science called fractal topology. Instead, we use the term fractal geometry, because the coordinates and analytic representations are usual ways of their construction. Mathematical morphology is a rather direct spinoff from topology to image processing. Using the hit or miss transformations we define dilation and erosion, opening and closing, etc. Again, it is the digital morphology (finite and discrete) that is used in computer algorithms for binary and gray images. The digital morphology had been extended to include random sets.

Most people associate topology with deterministic objects. While a number of papers were published about topology with fuzzy sets, the strength of fuzzy sets has been in control theory, not in modeling of timeless topological objects. While stochastic geometry had existed for many years, the stochastic topology -- the fusion of topological and stochastic methods -- is just in the beginning.

Clustering algorithms (including fuzzy clustering) belong to a very important branch of computational statistics called cluster analysis. The goal of cluster analysis is to find meaningful groups in multi-dimensional data, where each dimension (column of data) represents a feature of the measured population. Cluster analysis helps us to understand and sort out the real world data.

Shape descriptors of digital objects have been developed using various analytic techniques, for example, the Fourier transforms.

The graphics transformations known as morphing and warping are special cases of digitized homotopy mappings.

Gray geometry is a discrete geometry, where points have numerical values, for example, integer values from 0 to 255. Some algorithms calculate the boundary of a gray set. It can be, for example, the leading edge of a moving target. We need to be aware that in gray (discrete) geometry the boundary is ``thick'', not ``thin''. If one subtracts the boundary from a bounded set, then the diameter decreases. (Peeling an apple removes some mass and decreases its diameter.) Gray discrete geometry is a framework for digital image processing.

In 1979 Azriel Rosenfeld introduced a branch of Image Processing called Digital Topology, which is attached to classical topology at the conceptual level. This field is also known as mathematical morphology (Serra [1982, 1988], Najman \& Talbot [2010]). The major topic in digital topology is the study of connectedness, boundaries, edges, surfaces and manifolds. The book of Kong and Rosenfeld [1999] contains more than 360 references to papers in digital topology; however, none of these references is a paper in set-theoretic or general topology!

In digital topology, the discrete grid of pixels or voxels is regarded as a~sampling of a continuous image (a region in 2- or 3-dimensional Euclidean space, respectively). The pixels and voxels are at integer coordinates, so the grid is a lattice, and they have colors represented by numbers or vectors. From this viewpoint, it is of interest to study conditions under which this digitization process preserves topological (or other geometric) properties. This is important in manipulation of images created by Computer Tomography and Magnetic Resonance.

Image processing algorithms deal with digitized images. A digitized ima\-ge is the result of the analog to digital conversion; it is therefore a matrix of single or multiple values. The numbers in single-valued images usually represent intensities. If the values are 0 or 1, we have binary images. The gray scale images have many formats, for example, integer values from 0 to 255 (this is the 8-bit format). In sequence of images there is also a temporal aspect. The derivatives are actually differences, and these can be spatial or temporal. The images have features defined via point patterns and intensities. In feature extraction or pattern recognition we try to match an image to a subset of another image. In edge and corner detection, we try to identify and label these features in the image. In image segmentation we try to decompose an image into regions of interest. Target detection usually requires both pattern recognition and image segmentation. In image enhancement we try to improve visibility of certain patterns or objects. Image restoration requires some knowledge of noise or blur. Scene analysis uses logical reasoning schemes: many of the programs are examples of artificial intelligence programs.

The key to the applicability of a specific part of mathematics, for example, set-theoretic topology, to computer algorithms is a sufficient finite sampling of mathematical objects. We should keep in mind that for computer applications a finite sample of a topological space is not a discrete topological space. It is a representation of the space using a specified granularity related to a partition or to a covering. Using the computer graphics' language, it is a representation at certain level of detail. The points in finite samples are manipulated via their coordinates. The points in 2D representations are called pixels, and in 3D representations -- voxels.

There are two common misuses of the word ``topology''. The command ``tsort'' in UNIX/Linux system is called ``topological sorting'' (see also: Wikipedia), however, this sorting algorithm has nothing to do with topology. This algorithm takes a partial ordering relation on a finite set, and extends it to a linear ordering relation of the same set. Next, the common ``topologies'' in networks are called star, bus, ring, mesh, and cellular. Every mathematician would say that these are concepts of Graph Theory, not Topology.

\section{Fractal modeling}
Benoit Mandelbrot (1924--2010), who introduced the notion of fractals [1983], always pursued their real-world applications [1999]. Fractals are now used in physics, economics, electrical engineering, and financial mathematics.

Metric topology is only one of several sources of fractals. A few examples are the Cantor set and recursively constructed continua such as Sierpinski universal curves: the Sierpinski carpet and the Sierpinski gasket. The 3D version of Sierpinski carpet is called the Menger sponge. Recall also the space filling curves of Peano and Hilbert, and Koch curves. The impact of fractals on computer graphics is enormous. The notion of self-similarity qualifies the fractal objects to fractal geometry, not to fractal topology where the carpets, gaskets, etc. originated. Other important sources of fractals are the functions of complex variables (mappings and surfaces), and dynamical and chaotic systems.

While the symmetry of geometric objects is studied via groups of symmetry, the self-similarity is a somewhat more intrinsic notion. The affine mappings of an object into itself become part of its definition (construction). For example, the Cantor Set is covered by two smaller copies of the set (scaled 1:3). The finite approximations of fractals reflect their properties sufficiently to enable the computer applications.

The notion of fractal dimension (similar to Hausdorff dimension) is used in texture generation and texture analysis. Since this dimension can be calculated for parts of images, it is used for the classification of natural and synthetic textures.

The image compression algorithms, called fractal image compression, were developed and patented by Michael Barnsley and his collaborators (see [2006]). Several public domain fractal image compression algorithms (which are different from the patented algorithms) are available in the book of Yuval Fisher [1998], and at related web sites. The patented fractal image compression has been used in Microsoft Encyclopedia Encarta. The fractal compression leads to the notion of a fractal transform and the transform inverse. An advantage of the fractal compression technique is that the resized copies of an image (within a reasonable range) keep the qualities of the original size. So, in spite of loss of some data, the resized images appear much the same to the human eye.

Chaotic forms sometimes accompany the self-similarity in fractals. Hence we study the regular and chaotic behaviors in such dynamical systems. Fractal science has established the use of scientific visualization techniques in solution of iterative systems and differential equations. These are new avenues of the old branch of topology called topological dynamics.

Fractals are also an intellectual movement. The Internet contains web sites with truly beautiful and amazing images combining mathematical and artistic skills. Creating an image corresponding to a set of equations or functions just by iteration or recursion becomes the subject of a study by itself.

Iterative fractals can be calculated by parallel computers, while recursive (tree structured) fractals cannot. Some fractals are imploding (for example, the Sierpi\'nski gasket), and some are exploding (for example, the Pascal triangle). Some are approximating from the inside, like the Peano curve, and some from the outside, like the Cantor set. Some are deterministic like the Mandelbrot set, and some are random like the Barnsley fern. But ultimately, there is oneness of fractals in their fractal dimension.

Fractals enabled the stochastic modeling of difficult objects such as clouds, smokes, fires, plants, trees, terrain, and textures. However, in some cases, the results of a modeling should not be judged only by appearance. Better modeling requires deeper knowledge of physics, the measurement of parameters, etc. Some of these models study atmospheric turbulence and scintillation, which have demanding industrial and commercial applications.

Fractals tremendously accelerated the role of scientific visualization in the understanding of processes in nature described via mathematics. On the other side, we have images of amazing complexity created by imaging cameras: Earth from satellites, distant nebulas from the Hubble Space Telescope, and those from probes to extraterrestrial space.

\section{Moving forward}
Hugo Steinhaus (1887--1972), who founded Studia Mathematica with Stefan Banach, called for an establishment of a ``mathematical clinic'', where people would come for mathematical solutions. Perhaps, each math department should have at least one faculty member who connects mathematics to the external world (many already do).

It is inspiring to notice that early works of Waclaw Sierpinski (1882--1969) contain the ideas for recursive fractals [1916], and also for topological games [1924]. Please explore the web pages for his other important contributions.

Department of Geometry and Topology at Komensky University in Bratislava, Slovakia, was renamed to the Department of Computer Graphics and Image Processing. The new disciplines interact with the old ones. Other departments have also split resources between theoretical and practical work.

There are two kinds of computation: numerical and symbolic. Also, there are two kinds of descriptions of geometric objects: analytic and procedural. The analytic descriptions involve coordinates, functions and equations. The procedural descriptions are ``recipes'' which involve non-numerical algorithms, formal languages and shape grammars. For example, some fractals can be described using substitution rules with a simple alphabet.

Computer graphics is a practical branch of computational geometry, and has powerful applications in generating virtual reality scenes for the movies including interactive environments. Also, computational geometry algorithms are used in printers, plotters, and industrial robots. The software implementations of computational geometry are Computer Aided Design (CAD) and Computer Aided Manufacture (CAM). Other applications strongly influenced by computational geometry are the Geographic Information Systems (GIS). Recent applications of vector-valued images in remote sensing include multispectral and hyperspectral analysis. Statistics and especially multivariate analysis are important research tools in this area.

In truly applied mathematics, the criterion of beauty is not relevant. One has to be prepared for pages of ugly and long formulas, crude approximations, truncated constants, graphs of weird irregular curves, etc. A more useful mathematics is the one that yields to better modeling and forecasting, which are usually the domains of industrial physics and engineering.

The tracking of a moving object (passive or active) requires the measurement of the data, processing of the data, and controlling the devices. I~used a modified Hausdorff distance in software tracking multiple 3D objects. The mathematical framework for the consistent treatment of all these activities is known as the Kalman Filter, which is one of the most important algorithms of our time.

In computer modeling, the first and simple (na\"{\i}ve) algorithm usually does not work well. More mature algorithms result from multiple feedback cycles. Much as science uses experimentation to refine its theory to more faithfully reflect the real world.

The Ackermann function is computable, but not primitive recursive. However, it is practically non-computable, because we quickly run out of computer memory and computer time. You cannot turn the whole Earth into a~computer to get a few more values of the Ackermann function! The searches for decimals of $\pi$ and large prime numbers have also used a lot of computer memory and computer time. But, while in decimal system the decimal digits of $\pi$ are random, in hexadecimal number system we have recursive formula for its hexadecimal digits!

Formulas or equations are for most scientists and engineers a way to describe and understand a mathematical model. The model is built to explain some experimental data or to perform some desired actions. A different way is to use algorithms instead of formulas or equations. In particular, many researchers use software packages for interactive analysis and visualization of data. Cellular automata would be hard to develop, understand, and apply without such visualization.

As a mathematician (algorithm engineer) employed in industry I always asked myself some very basic questions. What is my task? What should I~accomplish? What is not my job, and when it is a task for someone else? What should I do earlier and what later? Who else knows a lot about it? It is important that I response to the immediate need of other scientists or engineers in the team. There is no time for making mathematical fictions! The funds are used for useful mathematical results.

The scientific work in an industrial environment depends on funding, a~degree of success in research and product development, and plans for future work. Funding does not come automatically, but must be obtained by knowing the demands, and by writing the proposals. A manager talks to a~client about work to be done. A manager talks to scientists and engineers how to approach the problem's solution. These experts take separate tasks, put things together, and the product is ready to be delivered. The funding scheme does not always support scientific logic. For example, the development of A was funded, the development of B was funded, but a comparison of A and B will not be funded, because there is neither time nor money to do it. Finding the opportunities for work requires entrepreneurial minds. See, for example, the web sites of Steven Wolfram (Mathematica) and Michael Barnsley (Iterated Systems) how they take advantage of the needs of the industrial community.

Most of mathematicians are spending their time trying to prove properties of objects that do not correspond to anything in the real world, while many phenomena of the real world remain unexplained, they are poorly modeled, or their features are not explored. We should be challenged by manufacture, agriculture, medicine and other activities aimed toward our survival and wellbeing.

\newpage

\begin{center}

{\bf Zastosowanie topologii w algorytmach komputerowych}
\end{center}
\medskip
\prefixing
\noindent{\bf Streszczenie.}
Celem niniejszego artyku/lu jest om\'owienie niekt/orych zastosowa\'n topologii og\'olnej w~algorytmach komputerowych w~kontek/scie modelowania i~symulacji, a~tak\.ze w~grafice komputerowej i~przetwarzaniu obraz/ow. Cho\'c post/ep w~tych dziedzinach jest bardzo uzale\.zniony od post/ep/ow w~sprz/ecie komputerowym, jednak g/l/ownymi osi/agni/eciami intelektualnymi s/a algorytmy. Zastosowania topologii og/olnej w~innych dzia/lach matematyki nie s/a omawiane, poniewa\.z nie s/a one zastosowaniami matematyki poza matematyk/a.\hfill\break
{\bf S/lowa kluczowe:} topologia og/olna, topologia teoriomnogo\'sciowa, topologia cyfrowa, matematyczna morfologia, przetwarzanie obraz/ow, modelowanie, symulacja, matematyka u\.zytkowa, matematyka stosowana.
\nonprefixing

\end{document}